\documentclass[oneside,leqno,11pt]{article}
\usepackage{revstat}
\usepackage{subcaption,graphicx,amsmath,epstopdf}%
\begin{document}
\title{On the Robustness and Asymptotic Properties for Maximum Likelihood Estimators of Parameters in Exponential Power and its Scale Mixture Form Distributions
\footnote{The opinions expressed in this text are those of the authors and do not necessarily reflect the views of any organization.}}
\renewcommand{\titleheading}
             {Using the ``revstat.sty" Package}  
\author{\authoraddress{Mehmet Niyazi \c{C}ankaya}
                      {Department of Statistics, Faculty of Arts and Science
                       U\c{s}ak University,\\
                       U\c{s}ak
                       \ (mehmet.cankaya@usak.edu.tr)}
                       \\
                        \authoraddress{Olcay Arslan}
                      {Department of Statistics, Faculty of  Science
                      Ankara University,\\
                      Ankara
                       \ (oarslan@ankara.edu.tr)}
      }
\renewcommand{\authorheading}
             {Author Abcdefg \ and \ Author Hijklmn}  

\maketitle

\begin{abstract}
The normality assumption on data set is very restrictive approach for modelling. The generalized form of normal distribution, named as an exponential power (EP) distribution, and its scale mixture form have been considered extensively to overcome the problem for modelling non-normal data set since last decades. However, examining the robustness properties of maximum likelihood (ML) estimators of parameters in these distributions, such as  the influence function, gross-error sensitivity, breakdown point and information-standardized sensitivity, has not been considered together. The well-known asymptotic properties of ML estimators of location, scale and added skewness parameters in EP and its scale mixture form distributions are studied and also these ML estimators for location, scale and scale variant (skewness) parameters can be represented as an iterative reweighting algorithm to compute the estimates of these parameters simultaneously. 
\end{abstract}

\begin{keywords}
Exponential power distributions; robustness; asymptotic; modelling.
\end{keywords}

\begin{ams}
49A05, 78B26.
\end{ams}


\section{Introduction}  
\label{Sec1}              


Modelling a data set from a real phenomena is an important issue to gain the best fitting on this data set. If there is a contamination in a data set, the inference procedure will be affected and the estimators will not be efficient and robust to the outlier(s) in the data set.  In this direction, the probability density functions (PDFs), especially that are an alternative to model the non-normally distributed data set, are overwhelmingly growing for the last decades. When there is a contamination in data set, a robust approach was considered by \cite{Hub64} to estimate the location parameter. After that time, important works were given by  \cite{Hampel68} and \cite{Hampel86} that give a comprehensive tools for testing the robustness properties of estimators.  The robustness properties of maximum likelihood (ML) estimators of parameters in the so-called PDFs should be proposed as well as an another approach in the framework of data modelling.  As a modelling issue, the asymmetry in data set was also considered by \cite{OHaganLeo76}. \cite{Azza85} proposed a method to have an asymmetric version of symmetric PDFs as a further processing on the modelling the asymmetry of data set if the asymmetry is in hand. The inference for the asymmetric normal distribution were studied by \cite{Azza85,Azza86,Henze86} and \cite{Chiogna98}. The asymmetric form of normal distribution in the sense of epsilon-skewed approach was proposed by \cite{Mud00}. Epsilon-skewed approach was considered to apply on the EP distribution by \cite{Elsal05}.  The bimodal case of epsilon-skew exponential power distribution for modelling the data set efficiently when the bimodality and asymmetry are observed together in data set was considered by \cite{cankaya15}. 


The second important problem is that the ML estimators of location, scale and skewness parameters in an arbitrary model are needed to compute. For this issue, the iterative reweighting algorithm (IRA) is applied by    ~\cite{Arslan04,Arslan95,kent91redescending,basu2004iteratively} and references therein to compute the estimates for location and scale parameters in an arbitrary model for the multivariate case.  The convergence performance of IRA to global point of the model has been approved by \cite{Arslan95,kent91redescending,Kent96,Arslan04,basu2004iteratively,ArslanGenc09} and reference therein. We will introduce the IRA as well. We will reconsider the IRA algorithm for the estimations of location, scale and skewness parameters. Thus, since the skewness parameter in our case is a scale variant form, it is seen that the IRA form can be constructed, as considered by \cite{ArslanGenc09} for the scale mixture form of the epsilon-skew exponential power (ESEP) distribution. The scale mixture form of ESEP generates a new function having a heavy-tailedness and asymmetry properties together, as proposed by \cite{ArslanGenc09}. We will examine the gross-error sensitivity (GES) that is defined to be an Euclidean norm of an influence function  (IF) as a local robustness, breakdown point (BDP) as a global robustness of ML estimate of location parameter and information-standardized sensitivity (ISS) of ML estimators of location, scale and skewness parameters in the ESEP and its scale mixture form distributions. Since the GES and ISS are functions of IF, they can be regarded as a local robustness property of estimator. One can read the papers about the modelling issue and the robustness properties of ML estimators for further details from Ref. \cite{Hampel86}. For example, a distribution proposed by \cite{ArslanGenc09} examining the robustness properties of ML estimators of parameters as well and \cite{gomezesgt} was considered to model the data set having heavy-tailedness and asymmetry together. The asymptotic properties of ML estimators of location, scale and skewness parameters from ESEP and its scale mixture form distributions are examined by means of the well-known regularity conditions. As a result, we were interested with distributions, an algorithm of computation for ML estimators, robustness of estimators and the asymptotic properties of ML estimators for the interested distributions that are ESEP and its scale mixture form distributions.



The sketch of paper is as follows. Section 2 represents the exponential power distribution and its scale mixture form briefly and gives ML estimators of location, scale and skewness parameters. Section 3 considers the robustness property of ML estimators, such as influence function, gross-error sensitivity,  breakdown point and information-standardized sensitivity. Section 4 includes the asymptotic properties of ML estimators of these parameters as a more detailed way. The section 5 is for the simulation and real data examples. The last section is for conclusion.


\section{Epsilon-Skew Exponential Power Distribution}\label{Sec2}

The PDF of Gamma is a generator to have PDF in an exponential power type, as proposed by \cite{Elsal05}. The epsilon-skew version of symmetric normal distribution as a special case of ESEP with $\alpha=2$ was considered by \cite{Mud00}. The PDF of ESEP is given as follow:
\begin{definition}
Let $X$ have $ESEP(\theta,\sigma,\varepsilon,\alpha)$ Then, $\theta \in \mathbb{R}$,  $\sigma >0$, $\varepsilon \in (-1,1)$ are location, scale and skewness parameters of model, respectively and  
  $\alpha>0$ is a parameter to change the peakedness of model. 
  \begin{equation}\label{generalizedepsilonskewnormal}
 f_{ESEP}(x)=\left\{
                \begin{array}{ll}
                  \frac{\alpha}{2\sqrt{2}\sigma\Gamma(1/\alpha)}  exp(-(\frac{\theta-x}{2^{1/2}(1+\varepsilon) \sigma})^\alpha),~ x < \theta  \; \\
                  \frac{\alpha}{2\sqrt{2}\sigma\Gamma(1/\alpha)}  exp(-(\frac{x-\theta}{2^{1/2}(1-\varepsilon) \sigma})^\alpha ),  ~ x \geq \theta .
                \end{array}
              \right.
\end{equation}
\end{definition}

If $X \sim  ESEP(\theta,\sigma,\varepsilon,\alpha)$, then  $r$th moment of random variable $X$ is
\begin{equation}\label{centralmomentESEP}
E(X - \theta)^r=\frac{2^{r/2}\sigma^r\Gamma((r+1)/\alpha)[(-1)^{-r}(1+\varepsilon)^{r+1}+(1-\varepsilon)^{r+1}]}{2\Gamma(1/\alpha)}, ~~~~ r >0 
\end{equation}
\noindent   $\int_0^\infty x^m e^{-\beta x^n}dx=\frac{\Gamma(\gamma)}{n\beta^\gamma}$, $\gamma=\frac{m+1}{n}$, $\beta>0, m>0, n>0$ is used to integrate  \cite{Gradshteyn07,Elsal05}.

\subsection{ML Estimators of Parameters $\theta$, $\sigma$  and $\varepsilon$ in ESEP Distribution}

Let $X_1, X_2, ...,X_n$ be independent and identical random variables and have $ESEP(\theta,\sigma,\varepsilon,\alpha)$. Here, the parameter $\alpha$ is fixed. The parameters $\theta$, $\sigma$  and $\varepsilon$ will be obtained by means of ML  estimation method as follow:
\begin{equation}\label{SNlikelihood}
\log(L(\theta,\sigma,\varepsilon;x))=n \log\bigg[\frac{\alpha}{2^{3/2}\Gamma(1/\alpha)\sigma}\bigg]-\sum_{i=1}^n \frac{|x_i -\theta|^\alpha}{[2^{1/2}(1-sign(x_i-\theta)\varepsilon)\sigma]^\alpha}
\end{equation}
\noindent The ML estimators will be obtained when log-likelihood function is maximized according to the parameters. The derivatives of log-likelihood function according to these parameters are taken, setting to zero and finally solving these equations according to the parameters simultaneously in order to get the ML estimators of parameters interested. 

\begin{eqnarray}\label{esepthetaest}
  \frac{\partial}{\partial \theta}\log L  &=& \sum_{i=1}^{n} \frac{\alpha|x_i - \theta|^{\alpha-1}sign(x_i - \theta)}{(\sqrt{2}(1-sign(x_i - \theta)\varepsilon)\sigma)^\alpha} = 0,  \\    
  \frac{\partial}{\partial \sigma}\log L&=& -\frac{n}{\sigma}+\sum_{i=1}^{n}\frac{\alpha|x_i - \theta|^\alpha}{(\sqrt{2}(1-sign(x_i - \theta)\varepsilon)\sigma)^\alpha\sigma}=0, \\  
  \frac{\partial}{\partial \varepsilon}\log L &=&\sum_{i=1}^{n}\frac{\alpha|x_i - \theta|^\alpha sign(x_i - \theta)}{(\sqrt{2}(1-sign(x_i - \theta)\varepsilon)\sigma)^\alpha(1-sign(x_i - \theta)\varepsilon)}=0.  
\end{eqnarray} 
The functions in equation \eqref{esepthetaest}-(2.6) are non-linear. We need to make a regularization on these equations to have a form that is appropriate to do the IRA algorithm for computing them. After algebraic manipulation is performed on these equations, the estimating equations will be obtained as follows:
\begin{eqnarray}\label{eq:wthetaESEPformula}
    \hat{\theta}  &=& \frac{\sum_{i=1}^n w(x_i) x_i}{\sum_{i=1}^n w(x_i)}, \\  
   \hat{\sigma}^{2}  &=&  \frac{1}{n}\sum_{i=1}^n w(x_i) (x_i - \hat{\theta})^2, \\  
    \hat{\varepsilon}  &=& \sum_{i=1}^n\frac{w(x_i)(x_i-\hat{\theta})^2sign(x_i-\hat{\theta})}{(1-sign(x_i-\hat{\theta})\hat{\varepsilon})^2}/\sum_{i=1}^n\frac{w(x_i)(x_i-\hat{\theta})^2}{(1-sign(x_i-\hat{\theta})\hat{\varepsilon})^2},
\end{eqnarray}
\noindent where the weight function $w$ is
\begin{flalign}\label{weightESEPformula}
~~~~~~~~~~~~~~ w(x_i)= \frac{\alpha |x_i - \hat{\theta}|^{\alpha-2}}{(2^{1/2}(1-sign(x_i - \hat{\theta})\hat{\varepsilon}))^\alpha \hat{\sigma}^{{\alpha-2}}}.
\end{flalign}
By using the equations \eqref{eq:wthetaESEPformula}-(2.9), the ML estimators of parameters $\theta$, $\sigma$  and  $\varepsilon$ are obtained \cite{phdcankaya15}. 

\subsection{Computation Steps of ML Estimators of the Parameters $\theta$, $\sigma$ and $\varepsilon$ in ESEP Distribution}

$X_n = \{x_1,x_2,...,x_n\}$ is a random sampling and $k \in \mathbb{N}^{+}$  is a number for representing iteration in IRA. Then, IRA steps are given by 

\textbf{1. Step} $\theta^{(1)}$, $\sigma^{(1)}$ and $\varepsilon^{(1)}$ are initial values of computation.

\textbf{2. Step} The value of weight function is computed as
 \begin{eqnarray*}
   w^{(k)}(x_i) &=& \frac{\alpha |x_i - \hat{\theta}^{(k)}|^{\alpha-2}}{(2^{1/2}(1-sign(x_i - \hat{\theta}^{(k)})\hat{\varepsilon}^{(k)}))^\alpha \hat{\sigma}^{{\alpha-2}^{(k)}}}  \nonumber
 \end{eqnarray*}

\textbf{3. Step}  The estimate of location parameter $\theta$ is
\begin{eqnarray*}
 \hat{\theta}^{(k+1)} &=& \frac{\sum_{i=1}^n w^{(k)}(x_i) x_i}{\sum_{i=1}^n w^{(k)}(x_i)}
\end{eqnarray*}

\textbf{4. Step} The estimate of scale parameter $\sigma$ is
\begin{eqnarray*}
\hat{\sigma}^{2^{(k+1)}} &=& \frac{1}{n}\sum_{i=1}^n w^{(k)}(x_i) (x_i - \hat{\theta}^{(k+1)})^2 
\end{eqnarray*}

\textbf{5. Step} The estimate of skewness parameter $\varepsilon$ is
\begin{eqnarray*}
  \hat{\varepsilon}^{(k+1)} &=& \sum_{i=1}^n\frac{w^{(k)}(x_i)(x_i-\hat{\theta}^{(k+1)})^2sign(x_i-\hat{\theta}^{(k+1)})}{(1-sign(x_i-\hat{\theta}^{(k+1)})\hat{\varepsilon}^{(k)})^2}/\sum_{i=1}^n\frac{w^{(k)}(x_i)(x_i-\hat{\theta}^{(k+1)})^2}{(1-sign(x_i-\hat{\theta}^{(k+1)})\hat{\varepsilon}^{(k+1)})^2}
  \end{eqnarray*}
\noindent Here, by using $\hat{\theta}^{(k+1)}$, $\hat{\sigma}^{(k+1)}$ and $\hat{\varepsilon}^{(k)}$,  the weight function in second step is recomputed as an updated weight in iteration.

\textbf{6. Step} If norm of vector $(\hat{\theta}^{(k+1)} - \hat{\theta}^{(k)},\hat{\sigma}^{(k+1)}-\hat{\sigma}^{(k)}, \hat{\varepsilon}^{(k+1)}-\hat{\varepsilon}^{(k)})^T$ is bigger than the prescribed value $e>0$, then the steps $2-5$ are repeated. Otherwise, the iteration is terminated and the last values of estimates in iteration are assigned to be estimates of these parameters \cite{phdcankaya15}.

\subsection{ESEP and its Scale Mixture Form Distributions}
The scale mixture form of ESEP is obtained by the variable transformation, considered by \cite{Theodo98,ArslanGenc09,gomezesgt,Sub23,Box62,Bala95book}. It is named as an epsilon-skew generalized t distribution (ESGT).

\begin{definition}
Let $X$ have $ESGT(\theta,\sigma,\varepsilon,\alpha,q)$. Then, $\theta \in \mathbb{R}$,  $\sigma >0$, $\varepsilon \in (-1,1)$ are location, scale and skewness parameters of model, respectively and  
  $\alpha>0$ and $q>0$  are parameters to change the peakedness and tail-thickness of model, respectively. 
\begin{equation}\label{epsilonSGT}
f_{ESGT}(x) = \frac{\alpha}{2\sqrt{2}B(1/\alpha,q)q^{1/\alpha}\sigma}\bigg(1+ \frac{|x - \theta|^\alpha}{2^{\alpha/2}(1-sign(x-\theta)\varepsilon)^\alpha q \sigma^\alpha} \bigg)^{-(\frac{\alpha q+1}{\alpha})}, x \in \mathbb{R}
\end{equation}
\end{definition}

If $X \sim  ESGT(\theta,\sigma,\varepsilon,\alpha)$, then  $r$th moment of random variable $X$ is
\begin{equation}\label{rthmomentofSGT}
   E((X-\theta)^r)=\frac{2^{r/2-1}q^{r/\alpha}\sigma^r\Gamma(\frac{r+1}{\alpha})\Gamma(q-r/\alpha)}{\Gamma(1/\alpha)\Gamma(q)}[(1-\varepsilon)^{r+1}+(1+\varepsilon)^{r+1}(-1)^{-r}], q\alpha>r
\end{equation}
\begin{equation}\label{expectedformula}
\int_0^\infty y^r (1+uy^k)^{-m} dy=\frac{B(\frac{r+1}{k},m-\frac{r+1}{k})}{ku^{\frac{r+1}{k}}}, ~~    0<\frac{r+1}{k}<m 
\end{equation}
\noindent is used to integrate \cite{Gradshteyn07}.

One can get the maximum likelihood estimators of the parameters $\theta,\sigma$ and $\varepsilon$. The ML estimators are same representation in the ML estimators of ESEP, so we omit to rewrite them. However, the weight function $w$ of ESt is as follow \cite{ArslanGenc09}:
\begin{equation}\label{west}
w(x_i) = \frac{\nu +1}{\nu (1-sign(x_i-\hat{\theta})\hat{\varepsilon})^2 + (\frac{x_i - \hat{\theta}}{\hat{\sigma}})^2}.
\end{equation} 

\section{Robustness Properties of ML Estimators of Parameters}
Before introducing the tools for robustness, the definition of M-estimator must be given to show that ML estimators are produced by their score functions. Then, for the ML estimators, the objective function is $\rho(x,\boldsymbol \tau)=-\log[f(x;\boldsymbol \tau)]$ and the score function is $\psi(x; \boldsymbol \tau) = \frac{\partial}{\partial \boldsymbol \tau} \rho(x; \boldsymbol \tau)$. The generalized version of ML estimation is a M-estimation and it is defined as follow \cite{Hampel86}:

\begin{definition}
A M-estimator is defined through a function $\rho:\mathcal{X} \times \Theta \rightarrow \mathbb{R}$ as the value $\hat{\boldsymbol \tau}(F) \in \mathbb{R}^p$ minimizing $\int \rho(x,\hat{\boldsymbol \tau})dF(x)$ over $\hat{\tau}$ or through a function $\psi:\mathcal{X} \times \Theta \rightarrow \mathbb{R}^p$ as the solution for $\hat{\boldsymbol \tau}$ of the vector equation $\int \psi(x_i, \hat{\boldsymbol \tau})dF(x)=0$.
\end{definition}

If the function $\psi$ has robustness properties, the ML estimators of parameters will be robust. The definition of influence function is as follow: 

\begin{equation}\label{IFofparametersvectors}
 IF(x;\hat{\boldsymbol \tau},F)=M(\psi,F)^{-1} \Psi(x;\hat{\boldsymbol \tau}).
\end{equation}
\noindent where $\boldsymbol \tau$ is a vector of parameters $\tau_1, \tau_2,...\tau_p$.  $\Psi = (\psi_1, \psi_2,...,\psi_p)$ is a vector of score functions derived by objective function $\rho$ after taking the derivatives with respect to (w.r.t) parameters of function $\rho$ and the matrix $M$ is
\begin{equation}\label{IFmultivariate}
  M(\psi,F)= -\int \bigg[\frac{\partial}{\partial \boldsymbol \tau} \psi(x; \boldsymbol   \tau) \bigg] f(x;\boldsymbol \tau)dx. 
\end{equation}
The robustness criteria generated by influence function are gross error and information-standardized sensitivities. The definitions of them are as follow \cite{Hampel74,Hampel86,Rieder94,JurecPic06,MarYo06,HubRoncbook09}:
\begin{definition}
The gross error sensitivity of estimators $\hat{\boldsymbol \tau}$ is
\begin{equation}\label{globalsens}
 \gamma_u^{*}(\hat{\boldsymbol \tau},F) =  \underset{x}{sup}  || IF(x;\hat{\boldsymbol \tau},F)||
\end{equation}
\end{definition}

\begin{definition}\label{selfInformationsensivitysection2}
   If $I(\tau)$ exists for all parameters in vector $\boldsymbol \tau$, then the information-standardized sensitivity is given by
\begin{equation}
\gamma_i^{*}(\hat{\boldsymbol \tau},F)   = \underset{x}{sup}  \{IF(x;\hat{\boldsymbol \tau},F)^TI( \boldsymbol \tau)IF(x;\hat{\boldsymbol \tau},F)\}^{1/2}.
\end{equation}
\end{definition}

The local robustness of ML estimators will be examined by the definition of influence function in equation \eqref{IFofparametersvectors}. The special cases of ESEP are epsilon-skew normal (ESN) for $\alpha=2$ and epsilon-skew Laplace (ESL) for $\alpha=1$. The following theorem shows whether or not the influence function ML estimators of parameters $\theta$, $\sigma$ and $\varepsilon$ will be bounded. It is noted that if the influence function is bounded, then ML estimators of these parameters will be local robust.

\subsection{Robustness Properties of ML Estimators of Parameters $\theta$, $\sigma$ and $\varepsilon$ in ESEP}

\begin{theorem}(Score functions)\label{IFESEP}
Let $X_1, X_2, ...,X_n$ be an independent and identically distributed random variables having $ESEP(\theta,\sigma,\varepsilon,\alpha)$. Then, the score functions of parameters  $\theta$, $\sigma$, $\varepsilon$ and $\alpha$ are as follows:
\begin{eqnarray}\label{ESEPpsis}
  \psi_\theta(x) &=& \frac{\alpha |x|^{\alpha-1}sign(x)}{[2^{1/2}(1-sign(x)\varepsilon)]^\alpha}, \\ \nonumber
  \psi_\sigma(x) &=& -1+\frac{\alpha|x|^\alpha}{[2^{1/2}(1-sign(x)\varepsilon)]^\alpha}, \\ \nonumber
   \psi_\varepsilon(x) &=&  \frac{\alpha|x|^{\alpha}sign(x)}{2^{\alpha/2}(1-sign(x)\varepsilon)^{\alpha+1}}, \\ \nonumber
 \psi_\alpha(x) &=&  -\frac{\Gamma(1/\alpha)+\Gamma^{'}(1/\alpha)\alpha^{-1}}{\alpha\Gamma(1/\alpha)}+\frac{|x|^\alpha\{\log|x|-\log[2^{1/2}(1-sign(x)\varepsilon)]\}}{[2^{1/2}(1-sign(x)\varepsilon)]^\alpha}. \nonumber
\end{eqnarray}
We will examine whether or not these score functions are finite when $x$ goes to infinity: 
\begin{enumerate}
\item[a)] For $\alpha >1$, $\underset{x \rightarrow \pm \infty}{\lim} \psi_\theta(x) = \infty$, $\underset{x \rightarrow \pm \infty}{\lim} \psi_\sigma(x) = \infty$, $\underset{x \rightarrow \pm \infty}{\lim} \psi_\varepsilon(x) = \infty$ and  $\underset{x \rightarrow \pm \infty}{\lim} \psi_\alpha(x) = \infty$. Then, the score functions of parameters $\theta$, $\sigma$, $\varepsilon$ and $\alpha$ are unbounded.
 \item[b)] For  $\alpha < 1$,  $\underset{x \rightarrow \pm \infty}{\lim} \psi_\theta(x) = 0$, $\underset{x \rightarrow \pm \infty}{\lim} \psi_\sigma(x) =  \infty$, $\underset{x \rightarrow \pm \infty}{\lim} \psi_\varepsilon(x) =  \infty$ and $\underset{x \rightarrow \pm \infty}{\lim} \psi_\alpha(x) =\infty$.
 \item[c)] For $\alpha = 1$,  $\underset{x \rightarrow \pm \infty}{\lim} \psi_\theta(x) = \frac{sign(x)}{2^{1/2}(1-sign(x)\varepsilon)}$,  $\underset{x \rightarrow \pm \infty}{\lim} \psi_\sigma(x) =  \infty$, $\underset{x \rightarrow \pm \infty}{\lim} \psi_\varepsilon(x) =  \infty$ and $\underset{x \rightarrow \pm \infty}{\lim} \psi_\alpha(x) =\infty$.
\end{enumerate}
\end{theorem}

\begin{corollary}(Influence Function)\label{commonIFESEP}
Let $X_1, X_2, ...,X_n$ be an independent and identically distributed random variables having $ESEP(\theta,\sigma,\varepsilon,\alpha)$. Then, the influence function of ML estimators for the parameters  $\theta$, $\sigma$ and $\varepsilon$ is
\begin{equation}\label{comIFESEPeqns}
IF(x;\hat{\theta},\hat{\sigma},\hat{\varepsilon}) =  \begin{bmatrix}
 T_{11}\psi_\theta(x)+ T_{12}\psi_\sigma(x)+ T_{13}\psi_\varepsilon(x) \\
 T_{21}\psi_\theta(x)+ T_{22}\psi_\sigma(x)+ T_{23}\psi_\varepsilon(x) \\
 T_{31}\psi_\theta(x)+ T_{32}\psi_\sigma(x)+ T_{33}\psi_\varepsilon(x)
  \end{bmatrix} = \begin{bmatrix} IF_{1}(x;\hat{\theta},\hat{\sigma},\hat{\varepsilon}) \\ IF_{2}(x;\hat{\theta},\hat{\sigma},\hat{\varepsilon}) \\ IF_{3}(x;\hat{\theta},\hat{\sigma},\hat{\varepsilon})   \end{bmatrix},
\end{equation}
\noindent where  $T_{ij}$ represents row $i$ and column $j$ of matrix $M$ in equation \eqref{IFmultivariate}, $i,j=1,2,3$. The ML estimators of these parameters are unbounded, because the score functions $\psi_\sigma(x)$ and $\psi_\varepsilon(x)$ are unbounded. As a result, the ML estimators of the parameters $\theta$, $\sigma$ and $\varepsilon$ are not robust.
\end{corollary}
Gross-Error Sensitivity: The norm of influence function of ML estimators for the parameters $\theta,\sigma$ and $\varepsilon$ of $ESEP(\theta,\sigma,\varepsilon)$  distribution is gross-error sensitivity.
\begin{equation}\label{GESESEPeqnsESEP}
 \gamma_u^*(\hat{\theta},\hat{\sigma},\hat{\varepsilon},F_{ESEP}) =\{(IF_{1})^2+(IF_{2})^2+(IF_{3})^2\}^{1/2}.
\end{equation}
Since the components $IF_1$, $IF_2$ or $IF_3$ of vector IF are not bounded, then  $\gamma_u^*(\hat{\theta},\hat{\sigma},\hat{\varepsilon},F_{ESEP})$ are not bounded.

The ISS in equation \eqref{selfInformationsensivitysection2} is obtained for the ML estimators of parameters $\theta$, $\sigma$ and $\varepsilon$ as follow:
\begin{eqnarray}\label{tseselfstandardized}
 \gamma_i^{*}(\hat{\theta},\hat{\sigma},\hat{\varepsilon}, F_{ESEP})  &=& \{  IF_1^2 I(\theta) + IF_2IF_1I(\sigma,\theta)+ IF_3 IF_1I(\varepsilon,\theta) \\  \nonumber
 &+& IF_1IF_2I(\theta,\sigma)  + IF_2^2I(\sigma)  + IF_3IF_2I(\varepsilon,\sigma) \\  \nonumber
  &+&    IF_1IF_3I(\theta,\varepsilon) + IF_2IF_3I(\sigma,\varepsilon) +  IF_3^2I(\varepsilon)\}^{1/2},
\end{eqnarray} 
\noindent where the components of  $IF_1,IF_2$ and $IF_3$ of vector $IF$ produce the functions $\psi_\theta^2, \psi_\theta\psi_\sigma,\psi_\theta\psi_\varepsilon,\psi_\sigma^2,\psi_\sigma\psi_\varepsilon$ and $\psi_\varepsilon^2$. If the results of these functions are bounded when $x \rightarrow \infty$ and also an each elements of Fisher matrix are bounded for the probable values of parameters, then ISS will be bounded.

\begin{corollary}(ISS)\label{propositionesepSelfInformation}
   Let  $X$ have  $ESEP(\theta,\sigma,\varepsilon,\alpha)$. The parameter $\alpha$ is fixed. Then,
    \begin{enumerate}
      \item For $\alpha=1$, $\underset{x \rightarrow  \infty}{\lim} \psi_\theta^2(x) = (2^{1/2}(1-sign(x)\varepsilon))^{-1}$.
      \item For $\alpha \in (0,1)$,  $\underset{x \rightarrow  \infty}{\lim} \psi_\theta^2(x)=0$.
      \item For $\alpha \in (1,\infty)$, $\underset{x \rightarrow  \infty}{\lim} \psi_\theta^2(x)=\infty$.
      \item For $\alpha=1/2$, $\underset{x \rightarrow  \infty}{\lim}  \psi_\theta(x) \psi_\varepsilon(x)=(4\sqrt{2}(1-sign(x)\varepsilon)^2)^{-1}$.
      \item For $\alpha \in (0,1/2)$, $\underset{x \rightarrow  \infty}{\lim} [\psi_\theta(x) \psi_\varepsilon(x)]=0$.
      \item For $\alpha \in (1/2,\infty)$, $\underset{x \rightarrow  \infty}{\lim} [\psi_\theta(x) \psi_\varepsilon(x)]=\infty$.
      \item When the functions $\psi_\theta^2(x)$ and $\psi_\theta(x)\psi_\varepsilon(x)$ in 2. and 5. are evaluated together, $\underset{x \rightarrow  \infty}{\lim} [\psi_\theta(x) \psi_\varepsilon(x)]=0$ and $\underset{x \rightarrow  \infty}{\lim} \psi_\theta^2(x)=0$  for $\alpha \in (0,1/2)$.
      \item When the functions $\psi_\theta(x) \psi_\varepsilon(x)$ and $\psi_\theta^2(x)$ in 3. and 6. are evaluated together,  $\underset{x \rightarrow  \infty}{\lim} [\psi_\theta(x) \psi_\varepsilon(x)]=\infty$, $\underset{x \rightarrow  \infty}{\lim} \psi_\theta^2(x)=\infty$  for  $\alpha \in (1,\infty)$.
     \item Since $\alpha>0$, $\underset{x \rightarrow  \infty}{\lim} \psi_\varepsilon^2(x)=\infty$, $\underset{x \rightarrow  \infty}{\lim} \psi_\sigma^2(x)=\infty$, $\underset{x \rightarrow  \infty}{\lim}  \psi_\sigma(x) \psi_\varepsilon(x) = \infty$.
     \item For $\alpha \in (0,1/2)$, $\underset{x \rightarrow  \infty}{\lim} \psi_\theta(x) \psi_\sigma(x)=0$, $\alpha =1/2$ için $\underset{x \rightarrow  \infty}{\lim} \psi_\theta(x) \psi_\sigma(x)=\frac{\alpha^2sign(x)}{(2^{1/2}(1-sign(x)\varepsilon))^{2\alpha}}$. 
     \item For $\alpha > 1/2$, $\underset{x \rightarrow  \infty}{\lim} \psi_\theta(x) \psi_\sigma(x)=\infty$.
     \item For $\alpha=1/2$, $\underset{x \rightarrow  \infty}{\lim} \psi_\theta(x)\psi_\varepsilon(x)<\infty$ . However, $\underset{x \rightarrow  \infty}{\lim}\psi_\varepsilon(x)=\infty$ for  $\alpha=1/2$.
  \end{enumerate}
 
\noindent As a result, the ML estimators of the parameters $\theta$, $\sigma$  and $\varepsilon$ of  $ESEP(\theta,\sigma,\varepsilon,\alpha)$ distribution are not robust in the sense of ISS. However, the score function of location parameter, namely $\underset{x \rightarrow  \infty}{\lim}  \psi_\theta(x)<\infty$, is finite when $x \rightarrow \infty$. Thus,  if the parameters $\sigma$ and $\varepsilon$ are fixed, the ML estimator of the parameter $\theta$ is robust.
\end{corollary} 

\begin{definition}\label{condiHubZhanLi}
Let $\psi$ be a score function. Suppose that there is a point $x_0$. For  $0<x \leq x_0$, $\psi$ is (weakly) increasing.
For $x_0 < x < \infty$, $\psi$ is (weakly) decreasing. Then, $\psi$ is defined to be a redescending function \cite{Huber84}.
\end{definition}

\begin{enumerate}
\item $\rho(0)=0$  \cite{Huber84,ZhangLi98}.
\item $\underset{|x| \rightarrow \infty}{\lim} \rho (x) = \infty$  \cite{Huber84,ZhangLi98}.
\item $\underset{|x| \rightarrow \infty}{\lim} \frac{\rho (x)}{|x|} = 0$ \cite{Huber84}.
\item  For  $0<x \leq x_0$, $\psi$ is weakly increasing. For $x_0 < x < \infty$, $\psi$ is weakly decreasing. There is a such a point $x_0$  \cite{Huber84}.
\end{enumerate}
The following condition is used to test whether or not function $\psi$ is monotone:
\begin{enumerate}
\item For  $0<x \leq x_0$, $\psi$ is a non decreasing function. For $x_0 < x < \infty$, $\psi$ is a non increasing function \cite{Huber84,ZhangLi98}.
\end{enumerate}
We will examine whether or not the conditions in definition \ref{condiHubZhanLi} will be satisfied for the ML estimator of location parameter $\theta$ in ESEP distribution. $\rho_{ESEP}=-\log(f_{ESEP})$ is objective function of ESEP distribution.

\begin{enumerate}
\item $\rho_{ESEP}(0)=0$.
\item $\underset{|x| \rightarrow \infty}{lim} \rho_{ ESEP}(x)=\infty$.
\item   $\underset{|x| \rightarrow \infty}{lim} \frac{\rho_{ ESEP}(x)}{|x|} =\left\{
                                                           \begin{array}{ll}
                                                             0, & \alpha<1; \\
                                                             \frac{1}{2^{1/2}(1-\varepsilon)}, & \alpha=1; \\
                                                             \infty, & \alpha>1.
                                                           \end{array}
                                                         \right.  $
\item In order to see that the function  $\psi_{ESEP}$  is non decreasing for the interval                                                   $0<x \leq x_0$ and non increasing for $x_0 < x < \infty$, $\psi_{ ESEP}^{'}(x)=\frac{\alpha(\alpha-1)|x|^{\alpha-2}}{(2^{1/2}(1-sign(x)\varepsilon))^\alpha}=0$ must be examined. The root of $\psi_{ ESEP}^{'}(x)$ according to $x$ is zero. This is obvious result of $\psi_{ ESEP}^{'}(x)$. Then, the condition $4$ is not satisfied. 
\end{enumerate}

Then, the following corollary for the breakdown point of ML estimator of $\theta$ is as follow:

\begin{corollary}(Breakdown Point)\label{BPESEP}
Let $X$ be distributed as  $ESEP(\theta,\sigma,\varepsilon,\alpha)$. For $\alpha>1$, since the conditions given above are not satisfied, the breakdown point of ML estimator of location parameter $\theta$ is not $1/2$, but for $\alpha \leq 1$, the breakdown point of ML estimator of $\theta$ is $1/2$. 
\end{corollary}

\subsection{Robustness Properties of ML Estimators of Parameters $\theta$, $\sigma$ and $\varepsilon$ in ESGT Distribution}

\begin{theorem}(Score functions)\label{IFESEP}
Let $X_1, X_2, ...,X_n$ be an independent and identically distributed random variables having  $ESGT(\theta,\sigma,\varepsilon,\alpha=2,q=\nu/2)$. Then, the score functions of parameters  $\theta$, $\sigma$, $\varepsilon$ and $\nu$ are as follows:
\begin{eqnarray}\label{Stpsis}
  \psi_\theta(x) &=& \frac{(\nu+1)x}{\nu[(1-sign(x)\varepsilon)]^2+x^2}, \\ \nonumber 
  \psi_\sigma(x) &=& \frac{(\nu+1) x^2}{\nu (1-sign(x)\varepsilon)^2+x^2} -1, \\ \nonumber 
  \psi_\varepsilon(x) &=& \frac{(\nu+1)x^2 sign(x)}{\nu(1-sign(x)\varepsilon)^3+(1-sign(x)\varepsilon) x^2}, \\ \nonumber
  \psi_\nu(x)  &=& k(\nu)+\frac{1}{2}\log\bigg[1+\frac{x^2}{\nu(1-sign(x)\varepsilon)^2}\bigg]-\frac{x^2}{\nu^2(1-sign(x)\varepsilon)^2+\nu x^2}\frac{\nu+1}{2}, \nonumber
\end{eqnarray}
\noindent where $k(\nu)=\frac{\Gamma^{'}(\frac{\nu+1}{2})0.5(\nu \pi)^{1/2}\Gamma(\nu/2)-\pi^{1/2}\Gamma(\frac{\nu+1}{2})[0.5\nu^{-1/2}\Gamma(\nu/2)+\Gamma^{'}(\nu/2)0.5\nu^{1/2}]}{(\nu \pi)^{1/2}\Gamma(\nu/2)\Gamma(\frac{\nu+1}{2})}$. In addition to, $ \underset{x \rightarrow  \pm \infty}{\lim} \psi_\theta(x) = 0$, $\underset{x \rightarrow  \pm \infty}{\lim} \psi_\sigma(x) = \nu$ and $\underset{x \rightarrow  \pm \infty}{\lim} \psi_\varepsilon(x) = ((\nu+1)sign(x))/(1-sign(x)\varepsilon)$. The score functions of parameters $\theta$, $\sigma$ and $\varepsilon$ are bounded, but  $\underset{x \rightarrow  \pm \infty}{\lim} \psi_\nu(x) = \infty$ \cite{ArslanGenc09}.
\end{theorem}

 Let $X_1, X_2, ...,X_n$ be an independent and identically distributed random variables having $ESGT(\theta,\sigma,\varepsilon,\alpha=2,q=\nu/2)$. Then, the influence function of ML estimators for the parameters  $\theta$, $\sigma$ and $\varepsilon$ is same with formula in equation \eqref{comIFESEPeqns}, as considered by \cite{ArslanGenc09}. The gross-error sensitivity in equation \eqref{GESESEPeqnsESEP} and information matrix (IM) in equation  \eqref{tseselfstandardized} are same for the ML estimators of parameters $\theta$, $\sigma$ and $\varepsilon$ in ESGT distribution. However, since the score functions of these parameters are finite, the IF of ML estimators of them is bounded. As it is defined that GES is a function of components $IF_1$, $IF_2$ and $IF_3$ of IF, the GES of ML estimators of them is bounded. Thus, ML estimators of them in ESGT will be robust if the parameters $\alpha$ and $q$ are taken to be fixed, as implied by \cite{ArslanGenc09}. 
 
\begin{corollary}(ISS)\label{propositionestSelfVarsens}
   Let  $X$ have  $ESGT(\theta,\sigma,\varepsilon,\alpha=2,q=\nu/2)$. The parameters $\alpha$ and $q$ are fixed. Then,
         \begin{enumerate}
       \item $\underset{x \rightarrow  \infty}{lim}\psi_\theta^2(x)=0.$  
       \item  $\underset{x \rightarrow  \infty}{lim}\psi_\theta(x)\psi_\sigma(x)=0.$  
       \item  $\underset{x \rightarrow  \infty}{lim}\psi_\theta(x)\psi_\varepsilon(x)=0.$  
       \item  $\underset{x \rightarrow  \infty}{lim}\psi_\sigma^2(x)=\nu^2.$  
       \item  $\underset{x \rightarrow  \infty}{lim}\psi_\sigma(x)\psi_\varepsilon(x)=\frac{\nu(\nu+1)}{1-sign(x)\varepsilon}.$  
       \item  $\underset{x \rightarrow  \infty}{lim}\psi_\varepsilon^2(x)=\frac{(\nu+1)^2}{(1-sign(x)\varepsilon)^2}.$  
     \end{enumerate}
     The IM given by \cite{gomezesgt} for ESGT distribution exists and each element of IM is finite when $\sigma<\infty,\nu<\infty$. Then, the ML estimators of parameters $\theta$, $\sigma$  and $\varepsilon$ are robust in the sense of ISS.
\end{corollary}

\subsubsection{The Breakdown point of ML estimator for location $\theta$ in $\varepsilon$-skew $t$ distribution}

The conditions given in definition \ref{condiHubZhanLi} are used to determine the breakdown point of ML estimator of location. Let $X$ be distributed as epsilon-skew $t$, abbreviated as $ ESt(\theta,\sigma,\varepsilon,\nu)$. $\rho_{ESt}=-\log(f)$ is objective function of ESt distribution.  The score function $\psi_{ESt}$ is derived from  $\rho_{ESt}=-\log(f_{ESt})$ after taking the derivative with respect to $\theta$. $\psi_{ESt}(x)=\frac{(\nu+1)x}{\nu(1-sign(x)\varepsilon)^2+x^2}$. 

\begin{enumerate}
\item $\rho_{ESt}(0)=0$,
\item $\underset{|x| \rightarrow \infty}{lim} \rho_{ESt}(x)=\infty$,
\item $\underset{|x| \rightarrow \infty}{lim} \frac{ \rho_{ESt}(x)}{|x|}=0,$
\item The derivative of the function $\psi_{ESt}(x)$ with respect to $x$ is $\psi_{ESt}^{'}(x)=\frac{(\nu+1)[\nu(1-sign(x)\varepsilon)^2+x^2]-2x^2(\nu+1)}{[\nu(1-sign(x)\varepsilon)^2+x^2]^2}$. The root of $\psi_{ESt}^{'}$ according to $x$ is $x_0 = \sqrt{\nu}(1-\varepsilon)$ for $x \geq 0$ and $x_0 = -\sqrt{\nu}(1+\varepsilon)$ for $x<0$. This function  increases on the interval $0 \leq x \leq \sqrt{\nu}(1-\varepsilon)$. It decreases on the interval  $-\sqrt{\nu}(1+\varepsilon) \leq x<0$. Then, we can find a root that function increases up to point $x_0$. After this point $x_0$, the function decreases. 
\end{enumerate}
Then, the following corollary for the breakdown point of ML estimator of location in ESt distribution is as follow:

\begin{corollary}(Breakdown Point)\label{BPEst}
Let $X$ be distributed as $ESt(\theta,\sigma,\varepsilon,\nu)$. Since the conditions given above are satisfied, the breakdown point of ML estimator of location parameter is $1/2$. 
\end{corollary}
 
\section{Asymptotic Properties of ML Estimators of the Parameters}
Let $X_1, X_2,...,X_n$  have a probability density function $f(x;\tau)$.  $\boldsymbol \tau =(\tau_1,\tau_2,...,\tau_p)^T$ is a vector of parameters and suppose the following regularity conditions hold \cite{Chanda54,LehmannCas98}:


\begin{description}
\item[(i)] $\boldsymbol\tau =(\tau_j,\tau_k,\tau_l)$, $j,k,l=1,2,...,p$. Then, $f(x;\boldsymbol\tau)$ is defined and $\tau_j,\tau_k,\tau_l$ are completely specified. Otherwise, the estimators of parameters are not consistent, i.e., if  $\boldsymbol \tau \neq \boldsymbol \tau^{'}$, then $f(x; \boldsymbol \tau) \neq f(x; \boldsymbol \tau^{'})$.
\item[(ii)] The distributions $f(x; \boldsymbol \tau)$ have common support.
\item[(iii)] The observations are $X_n=\{x_1,x_2,...,x_n\}$, where the $x_i$ are (independent identical distribution) iid with probability density $f(x; \boldsymbol \tau)$ with respect to $\mu$.
     \item[(iv)] The parameter space $\Omega$ contains an open set $\omega$  of which the true parameter value $\tau_0$ is an interior point. For each  $\tau_j,\tau_k,\tau_l \in \omega$, $\frac{\partial }{\partial \tau_j}\log f(x;\tau_j,\tau_k,\tau_l)$,  $\frac{\partial^2 }{\partial \tau_i \partial \tau_j}\log f(x;\tau_j,\tau_k,\tau_l)$ and $\frac{\partial^3}{\partial \tau_j \partial \tau_k \partial \tau_l} \log f(x;\tau_j,\tau_k,\tau_l)$ exist.
    \item[(v)]  $E[\frac{\partial}{\partial \boldsymbol \tau} \log f(x; \boldsymbol \tau)]=0$, $I(\boldsymbol \tau) = E[-\frac{\partial^2}{\partial\tau_j \partial \tau_k} \log f(x;\boldsymbol \tau)]$. $I(\boldsymbol \tau)$  is Fisher information matrix.
  \item[(vi) ]  $det(I(\boldsymbol \tau)) <\infty$ is satisfied.
  \item[(vii) ] Suppose that there exist functions $M_{jkl}$ such that
    $$|\frac{\partial^3}{\partial \tau_j \partial \tau_k \partial \tau_l} \log f(x;\boldsymbol \tau)| \leq M_{jkl}(x).$$ 
    $E[M_{jkl}(X)]< \infty$ is satisfied.
\end{description}

\begin{theorem}\label{LehmannCasellatheorem5.1p463}
  Let $X_1,...,X_n$ be iid, each with a density $f(x;\boldsymbol \tau)$ (with respect to
$\mu$) which satisfies the conditions above. Then, with probability tending to 1 as $n \rightarrow \infty$, there exist solutions $\hat{\boldsymbol \tau}$ of the log-likelihood equations derived from $\log(f)$ after taking derivative w.r.t parameters such that
    \begin{enumerate}
  \item $\hat{\tau}_{j}$ is consistent for estimating $\tau_j$,
 \item  $\sqrt{n}(\hat{\boldsymbol \tau}- \boldsymbol \tau)$      is asymptotically normal with (vector) mean zero and covariance matrix $[I(\mathbf{\boldsymbol \tau})]^{-1}$, and
  \item $\hat{\tau}_{j}$ is asymptotically efficient and asymptotic normally distributed.
    \end{enumerate}
\end{theorem}

\subsection{Asymptotic Properties of ML Estimators of the Parameters in ESEP Distribution}

Let $X$ have $ESEP(\theta,\sigma,\varepsilon,\alpha)$. The parameter $\alpha$ is fixed. The Fisher information matrix for the parameters $\theta$, $\sigma$ and $\varepsilon$ is

\begin{flalign*}
 \mathbf{I}_{ESEP}  = n
 \begin{bmatrix}
 \frac{\alpha(\alpha-1)\Gamma(1-1/\alpha)}{2\sigma^2\Gamma(1/\alpha)(1-\varepsilon^2)}& 0  & \frac{\alpha^2}{\sqrt{2}\sigma\Gamma(1/\alpha)(1-\varepsilon^2)} \\
   & \frac{-1}{\sigma^2}+\frac{\alpha(\alpha+1)\Gamma(1+1/\alpha)}{\Gamma(1/\alpha)\sigma^2} & 0 \\
   &  &   \frac{\alpha(\alpha+1)\Gamma(1+1/\alpha)}{\Gamma(1/\alpha)(1-\varepsilon^2)}
 \end{bmatrix}.
\end{flalign*}

The asymptotic variance-covariances of ML estimators of the parameters $\theta$, $\sigma$  and $\varepsilon$ are as follows:
\begin{small}
\begin{eqnarray}
  Var(\hat{\theta}) &=&  \frac{2(\alpha+1)\Gamma(1+1/\alpha)\Gamma(1/\alpha)\sigma^2(1-\varepsilon^2)}{\alpha n[\Gamma(1+1/\alpha)\Gamma(1-1/\alpha)\alpha^2-\Gamma(1-1/\alpha)\Gamma(1+1/\alpha)-\alpha^2]}, \\ \nonumber
\\
   Cov(\hat{\theta},\hat{\sigma}) &=& 0, \\  \nonumber
 \\
 Cov(\hat{\theta},\hat{\varepsilon}) &=&  \frac{\sqrt{2}\sigma\Gamma(1/\alpha)(\varepsilon^2-1)}{n[\Gamma(1+1/\alpha)\Gamma(1-1/\alpha)\alpha^2-\Gamma(1-1/\alpha)\Gamma(1+1/\alpha)-\alpha^2]},\\ \nonumber
 \\
  Var(\hat{\sigma})  &=&-\frac{\Gamma(1/\alpha)\sigma^2}{n[-\alpha^2\Gamma(1+1/\alpha)-\alpha\Gamma(1+1/\alpha)+\Gamma(1/\alpha)]}, \\ \nonumber
  \\
 Cov(\hat{\sigma},\hat{\varepsilon})  &=&0, \\ \nonumber
 \\
Var(\hat{\varepsilon})  &=& \frac{(\alpha-1)\Gamma(1-1/\alpha)\Gamma(1/\alpha)(1-\varepsilon^2)}{\alpha n[\Gamma(1+1/\alpha)\Gamma(1-1/\alpha)\alpha^2-\Gamma(1-1/\alpha)\Gamma(1+1/\alpha)-\alpha^2] }.
\end{eqnarray}
\end{small}
The Fisher information matrix is calculated again for the ESEP distribution. The variance-covariance matrix in Ref. \cite{Mud00} is same with our inverse of Fisher information matrix of ESEP with $\alpha=2$.

 \begin{table}[htb]
 \label{ESNvarcovvalues}
 \centering
 \caption{The asymptotic variance values of ML estimators for parameters $\theta,\sigma$ and $\varepsilon$ in ESN distribution with $\alpha=2$ }
\begin{tabular}{cccccc}
\hline
                   & $Var(\hat{\tau})/n$ &                             $n=30$ & $n=50$ & $n=100$ & $n=150$ \\ \hline \hline
                   & $Var(\hat{\theta})/n$                            & $0.211680 $    & $0.127008 $    & $0.063504 $     & $0.042336 $     \\
$\varepsilon=-0.2$ & $Var(\hat{\sigma})/n$                 & $0.016667 $    & $0.010000 $    & $0.005000 $     & $0.003333 $     \\
                   & $Var(\hat{\varepsilon})/n$                     & $0.070560 $    & $0.042336 $    & $0.021168 $     & $0.014112 $     \\ \hline

                   & $Var(\hat{\theta})/n$                            & $0.165375 $    & $0.099225 $    & $0.049612 $     & $0.033075 $     \\
$\varepsilon=-0.5$ & $Var(\hat{\sigma})/n$                 & $0.016667 $    & $0.010000 $    & $0.005000 $     & $0.003333 $     \\
                   & $Var(\hat{\varepsilon})/n$                     & $0.055125 $    & $0.033075 $    & $0.016538 $     & $0.011025 $     \\ \hline

                   & $Var(\hat{\theta})/n$                            & $0.079380 $    & $0.047628 $    & $0.023814 $     & $0.015876 $     \\
$\varepsilon=-0.8$ & $Var(\hat{\sigma})/n$                 & $0.016667 $    & $0.010000 $    & $0.005000 $     & $0.003333 $     \\
                   & $Var(\hat{\varepsilon})/n$                     & $0.026460 $    & $0.015876 $    & $0.007938 $     & $0.005292 $     \\ \hline
\end{tabular}
\end{table}



When the conditions $(iv)$, $(v)$, $(vi)$ and $(vii)$ are satisfied, the ML estimators of the parameters $\theta$, $\sigma$ and $\varepsilon$ for ESN distributions are asymptotic normally distributed, consistent and asymptotic efficient. In order to see whether or not the condition $(iv)$ is satisfied, Mathematica codes in appendix can be used. One observes that the derivatives of parameters of ESEP can be obtained. The condition $v$ showing that Fisher information matrix exist is satisfied. The condition $(vi)$ is satisfied, because 
$\det(\mathbf{I}_{ESEP})=\frac{(\Gamma(1-1/\alpha)\Gamma(1+1/\alpha)\alpha^2-\Gamma(1-1/\alpha)\Gamma(1+1/\alpha)\alpha^2-\Gamma(1-1/\alpha)\Gamma(1+1/\alpha)-\alpha^2)}{n^{-3}\alpha^{-2}(-2)\sigma^4\Gamma(1/\alpha)^3(1-\varepsilon^2)^2(-\alpha^2\Gamma(1+1/\alpha)-\alpha\Gamma(1+1/\alpha)+\Gamma(1/\alpha))^{-1}}<\infty.$
For $\alpha=2$, one can get $det(\mathbf{I}_{ESN})=\frac{2n^3(3\pi-8)}{\sigma^4\pi(1-\varepsilon^2)^2}<\infty$. The condition $(vii)$ is satisfied, because one can get a function $M_{jkl}=(X- \theta)^r$. From equation  \eqref{centralmomentESEP}, $E(X-\theta)^r$ exists and is finite.

The special case of ESEP with $\alpha=1$ is Laplace distribution. One can read that the PDF of Laplace distribution does not satisfy the regularity conditions \cite{Hoggetal}.  It is noted that since the epsilon-skew form is a kind of scaling, the property of Laplace distribution does not change. Due to this reason, the asymptotic properties of epsilon-skew Laplace (ESL) is not examined by means of classical derivative. The same situation is for $\alpha<1$. In future, we will examine this function by means of tools in fractional calculus \cite{Baletal12}.

\subsection{Asymptotic Properties of ML Estimators of the Parameters in ESGT Distribution}

The scale mixture form of ESEP distribution is studied by \cite{ArslanGenc09,gomezesgt}.  We got this distribution again. When our ESGT and the ESGT in \cite{gomezesgt} are compared, it is observed that they are same, because in our ESGT, there is a $\sqrt{2}$ part. It can be considered as a rewritten form of parameter q. So, instead of using PDF of our ESGT, we will use the ESGT of \cite{ArslanGenc09,gomezesgt}. It is noted that we prefer to keep the property of ESEP (because the normal distribution can be dropped for the special values of parameters in ESEP such that $\alpha=2$ and $\varepsilon=0$), so we proposed our ESGT in equation \eqref{epsilonSGT} again to show that it is easily obtained from the scale mixture form of distribution in ESEP type when $\nu$ goes to infinity, as a well-known property between normal and $t$ distributions. 

In this work, we give an additional step.  The Fisher information matrix of ESGT distribution was obtained by \cite{gomezesgt}. However, one integral in Fisher information was not be calculated. The result of that integral is $\int_0^\infty (1+y^\alpha)^{(-q-1/\alpha)} \log(1+y^p)dy=-\frac{\Gamma(1/\alpha)\Gamma(q)(\psi(q) - \psi(q+1/\alpha))}{\Gamma(q+1/\alpha)\alpha}$. To calculate this integral, one can take the derivative of the integral formula in equation \eqref{expectedformula} with respect to $m$.  Thus, the finiteness of determinant (condition $(vi)$) of Fisher information for ESGT distribution is guaranteed (see section 4 in \cite{gomezesgt}). The determinant of Fisher information for ESGT with $\alpha=2$ and $q=\nu/2$ distribution is $-1/2\frac{n^3(\nu+1)^2\nu(16 c^2-3)}{\sigma^6 (\varepsilon^2-1)^2(\nu+3)^3}$, where $c=\frac{\Gamma(\frac{\nu+1}{2})}{\sqrt{\nu \pi}\Gamma(\nu/2)}$ is a normalizing constant for the condition $(vi)$. The condition $(iv)$ can be satisfied from the Mathematica codes in appendix. To satisfy the condition $(vii)$, $M_{jkl}=(X- \theta)^r$ can be taken. From equation  \eqref{rthmomentofSGT}, $E(X-\theta)^r$ exists and is finite.

The first part of condition $v$, i.e, $E[\frac{\partial}{\partial \boldsymbol \tau} \log f(x; \boldsymbol \tau)]=0$, is satisfied for ESEP and its scale mixture form, because the support of the functions ESEP and ESGT does not depend on the parameter, as known from the regularity conditions \cite{Hoggetal}. 

\begin{table}[htb]
\label{EStvarcovvalues}
\centering
\caption{The asymptotic variance values of ML estimators for parameters $\theta,\sigma$ and $\varepsilon$ in ESt distribution with $\nu=3$ }
\begin{tabular}{cccccc}
\hline
                   & $Var(\hat{\tau})/n$ & $n=30$ & $n=50$ & $n=100$ & $n=150$ \\ \hline \hline
                   & $Var(\hat{\theta})/n$                 & $0.17173 $    & $0.103040 $    & $0.051520 $     & $0.034347 $     \\
$\varepsilon=-0.2$ & $Var(\hat{\sigma})/n$                 & $0.13333 $    & $0.080000 $    & $0.040000 $     & $0.026667 $     \\
                   & $Var(\hat{\varepsilon})/n$            & $0.05725 $    & $0.034347 $    & $0.017173 $     & $0.011449 $     \\ \hline
                   & $Var(\hat{\theta})/n$                 & $0.13417 $    & $0.080500 $    & $0.040250 $     & $0.026834 $     \\
$\varepsilon=-0.5$ & $Var(\hat{\sigma})/n$                 & $0.13333 $    & $0.080000 $    & $0.040000 $     & $0.026667 $     \\
                   & $Var(\hat{\varepsilon})/n$            & $0.04472 $    & $0.026834 $    & $0.013417 $     & $0.008945 $     \\ \hline
                   & $Var(\hat{\theta})/n$                 & $0.06440 $    & $0.038640 $    & $0.019320 $     & $0.012880 $     \\
$\varepsilon=-0.8$ & $Var(\hat{\sigma})/n$                 & $0.13333 $    & $0.080000 $    & $0.040000 $     & $ 0.026667$     \\
                   & $Var(\hat{\varepsilon})/n$            & $0.02147 $    & $0.012880 $    & $0.006440 $     & $ 0.004293$     \\ \hline
\end{tabular}
\end{table}

\section{Simulation and Real Data}

The simulation is considered to see the performance of the asymmetry parameter in model and its simultaneous estimation with location and scale parameters via using the IRA. We also provide the real data examples as an illustration of modelling capacity when there is an asymmetry in data set. The goodness of fit test Kolmogorov-Smirnov (KS), Akaike and Bayesian information criterions, abbreviated as AIC and BIC, respectively, are given to see the fitting performance of the ESN, ESL and ESt distributions.

Three cases are considered in the simulation study. In these three cases, the degrees of asymmetry are taken as low ($\varepsilon=-0.2$), middle ($\varepsilon=-0.5$) and high ($\varepsilon=-0.8$).  The random numbers are generated from its functions, because the performance of ML estimators for four sample sizes $n=30,50,100$ and $150$ from ESN, ESL and ESt distributions when the parameters $\theta$, $\sigma$ and $\varepsilon$ are estimated simultaneously is considered to test. To see the performance of simultaneous estimation via IRA, one criterion that is mean squared error (MSE) from simulation is used. The variance of each estimator from simulation is also given. The number of replication is $1000$.

In the simulation, there are three cases as follows:

\begin{itemize}
\item \textbf{Case I}: $ESN(\theta=0,\sigma=1,\varepsilon=\varepsilon_0)$
\item \textbf{Case II}: $ESL(\theta=0,\sigma=1,\varepsilon=\varepsilon_0)$
\item \textbf{Case III}: $ESt(\theta=0,\sigma=1,\varepsilon=\varepsilon_0)$
\end{itemize}

The algorithm based on the variable transformation is used to generate the random numbers from corresponding distributions. 

To generate the random numbers from ESN and ESL distributions, the following algorithm is used in order:
\begin{enumerate}
\item The random numbers are generated from $g \sim \Gamma(1/\alpha,1)$ function.
\item The random numbers are generated from $u \sim U(0,1)$ function.
\item If $u< \frac{1-\varepsilon}{2}$, then $ i=\sqrt{2}(1-\varepsilon)$, else $i=-\sqrt{2}(1+\varepsilon)$.
\item Using $x=\theta + \sigma i g^{1/\alpha}$, the random numbers are generated from ESN and ESL functions for $\alpha=2$ and  $\alpha=1$, respectively.
\end{enumerate}
To generate the random numbers from ESt distribution, the following algorithm is used in order:
\begin{enumerate}
\item The random numbers are generated from $g \sim \Gamma(1/2,1)$ function.
\item The random numbers are generated from $u \sim U(0,1)$ function.
\item If $d< \frac{1-\varepsilon}{2}$, then $ i=\sqrt{2}(1-\varepsilon)$, else $i=-\sqrt{2}(1+\varepsilon)$.
\item Using $y=\theta + \sigma i g^{1/2}$, the random numbers are generated from ESN distribution. 
\item The random numbers are generated from $z \sim \Gamma(\nu/2,1)$, where $\nu=3$.
\item Using  $x = y z^{-1/2} (\nu/2)^{1/2}$, the random numbers are generated. Thus, the random numbers are generated from ESt function.
\end{enumerate}

\begin{table}[!htb]
\centering
\caption{Case I: The ML estimates of parameters $\theta$,  $\sigma$ and $\varepsilon$ from ESN distribution for $\varepsilon=-0.2, -0.5$ and $-0.8$.}
\label{simESN}
\scalebox{0.9}{
\begin{tabular}{cccccccc}\hline
& $\tau$ & $\hat{\tau}$ & ${Var}(\hat{\tau})$ & ${MSE}(\hat{\tau})$  & $\hat{\tau}$ & ${Var}(\hat{\tau})$ & ${MSE}(\hat{\tau})$  \\  \hline \hline
 &   & & $n=      30$ & &   &$n=      50$  &   \\ \hline\hline
$\theta$&$0.0$&$ -0.0055$&$  0.1411$&$  0.1411$&$ -0.0250$&$  0.0765$&$  0.0771$\\

$\sigma$&$1.0$&$  0.9535$&$  0.0161$&$  0.0182$&$  0.9721$&$  0.0107$&$  0.0115$\\

$\varepsilon$&$-0.2$&$ -0.2185$&$  0.0496$&$  0.0499$&$ -0.2150$&$  0.0298$&$  0.0300$\\
\hline\hline
 &   & & $n=     100$ & &  &$n=     150$  &   \\ \hline\hline
$\theta$&$0.0$&$  0.0531$&$  0.0421$&$  0.0449$&$  0.0008$&$  0.0249$&$  0.0249$\\

$\sigma$&$1.0$&$  0.9846$&$  0.0048$&$  0.0050$&$  0.9980$&$  0.0033$&$  0.0033$\\

$\varepsilon$&$-0.2$&$ -0.1786$&$  0.0139$&$  0.0143$&$ -0.2050$&$  0.0087$&$  0.0088$\\
\hline \hline
 &   & & $n=      30$ & &   &$n=      50$  &   \\ \hline\hline
 $\theta$&$0.0$&$  0.0211$&$  0.1140$&$  0.1145$&$ -0.0061$&$  0.0681$&$  0.0682$\\

$\sigma$&$1.0$&$  0.9703$&$  0.0175$&$  0.0184$&$  0.9772$&$  0.0101$&$  0.0106$\\

$\varepsilon$&$-0.5$&$ -0.5040$&$  0.0373$&$  0.0373$&$ -0.5108$&$  0.0234$&$  0.0235$\\ \hline \hline
 &   & & $n=     100$ & &  &$n=     150$  &   \\ \hline\hline
$\theta$&$0.0$&$  0.1330$&$  0.0328$&$  0.0504$&$  0.1131$&$  0.0204$&$  0.0332$\\

$\sigma$&$1.0$&$  1.0001$&$  0.0050$&$  0.0050$&$  1.0045$&$  0.0033$&$  0.0034$\\

$\varepsilon$&$-0.5$&$ -0.5214$&$  0.0225$&$  0.0229$&$ -0.5127$&$  0.0061$&$  0.0081$\\
\hline\hline
 &   & & $n=      30$ & &   &$n=      50$  &   \\ \hline\hline
$\theta$&$0.0$&$  0.1721$&$  0.0846$&$  0.2231$&$  0.1337$&$  0.0599$&$  0.1712$\\

$\sigma$&$1.0$&$  1.0188$&$  0.0215$&$  0.0218$&$  1.0215$&$  0.0102$&$  0.0107$\\

$\varepsilon$&$-0.8$&$ -0.7223$&$  0.0415$&$  0.0476$&$ -0.7487$&$  0.0274$&$  0.0301$\\ \hline\hline
 &   & & $n=     100$ & &  &$n=     150$  &   \\ \hline\hline
$\theta$&$0.0$&$  0.0233$&$  0.0252$&$  0.0357$&$  0.0183$&$  0.0102$&$  0.0103$\\

$\sigma$&$1.0$&$  1.0222$&$  0.0062$&$  0.0066$&$  1.0243$&$  0.0038$&$  0.0043$\\

$\varepsilon$&$-0.8$&$ -0.7530$&$  0.0177$&$  0.0199$&$ -0.7716$&$  0.0076$&$  0.0084$\\
\hline\hline
\end{tabular}
}
\end{table}
\begin{table}[!htb]
\centering
\caption{Case II: The ML estimates of parameters $\theta$,  $\sigma$ and $\varepsilon$ from ESL distribution for $\varepsilon=-0.2, -0.5$ and $-0.8$.}
\label{simESL}
\scalebox{0.9}{
\begin{tabular}{cccccccc}\hline
& $\tau$ & $\hat{\tau}$ & ${Var}(\hat{\tau})$ & ${MSE}(\hat{\tau})$  & $\hat{\tau}$ & ${Var}(\hat{\tau})$ & ${MSE}(\hat{\tau})$  \\  \hline \hline
 &   & & $n=      30$ & &   &$n=      50$  &   \\ \hline\hline
$\theta$&$0.0$&$  0.0443$&$  0.0416$&$  0.0435$&$  0.0396$&$  0.0203$&$  0.0219$\\

$\sigma$&$1.0$&$  0.9821$&$  0.0297$&$  0.0300$&$  0.9731$&$  0.0196$&$  0.0203$\\

$\varepsilon$&$-0.2$&$ -0.1964$&$  0.0190$&$  0.0190$&$ -0.1888$&$  0.0125$&$  0.0126$\\
\hline\hline
 &   & & $n=     100$ & &  &$n=     150$  &   \\ \hline\hline
$\theta$&$0.0$&$  0.0395$&$  0.0095$&$  0.0110$&$  0.0319$&$  0.0055$&$  0.0065$\\

$\sigma$&$1.0$&$  0.9971$&$  0.0100$&$  0.0100$&$  0.9936$&$  0.0065$&$  0.0065$\\

$\varepsilon$&$-0.2$&$ -0.1845$&$  0.0051$&$  0.0053$&$ -0.1888$&$  0.0035$&$  0.0036$\\
\hline \hline
 &   & & $n=      30$ & &   &$n=      50$  &   \\ \hline\hline
$\theta$&$0.0$&$  0.1163$&$  0.0314$&$  0.0449$&$  0.0774$&$  0.0182$&$  0.0242$\\

$\sigma$&$1.0$&$  0.9635$&$  0.0235$&$  0.0248$&$  1.0045$&$  0.0235$&$  0.0235$\\

$\varepsilon$&$-0.5$&$ -0.4675$&$  0.0158$&$  0.0169$&$ -0.4564$&$  0.0076$&$  0.0095$\\ 
 \hline \hline
 &   & & $n=     100$ & &  &$n=     150$  &   \\ \hline\hline
$\theta$&$0.0$&$  0.0810$&$  0.0112$&$  0.0177$&$  0.0779$&$  0.0060$&$  0.0121$\\

$\sigma$&$1.0$&$  0.9960$&$  0.0090$&$  0.0090$&$  0.9922$&$  0.0077$&$  0.0078$\\

$\varepsilon$&$-0.5$&$ -0.4691$&$  0.0035$&$  0.0044$&$ -0.4675$&$  0.0023$&$  0.0033$\\
\hline\hline
 &   & & $n=      30$ & &   &$n=      50$  &   \\ \hline\hline
 $\theta$&$0.0$&$  0.1810$&$  0.0466$&$  0.0794$&$  0.1757$&$  0.0160$&$  0.0469$\\

$\sigma$&$1.0$&$  0.9392$&$  0.0222$&$  0.0259$&$  0.9993$&$  0.0155$&$  0.0155$\\

$\varepsilon$&$-0.8$&$ -0.7598$&$  0.0050$&$  0.0066$&$ -0.7456$&$  0.0050$&$  0.0060$\\
 \hline\hline
 &   & & $n=     100$ & &  &$n=     150$  &   \\ \hline\hline
$\theta$&$0.0$&$  0.1692$&$  0.0075$&$  0.0361$&$  0.1575$&$  0.0095$&$  0.0343$\\

$\sigma$&$1.0$&$  0.9686$&$  0.0102$&$  0.0109$&$  1.0070$&$  0.0056$&$  0.0056$\\

$\varepsilon$&$-0.8$&$ -0.7380$&$  0.0046$&$  0.0059$&$ -0.7494$&$  0.0021$&$  0.0046$\\
\hline\hline
\end{tabular}
}
\end{table}
\begin{table}[!htb]
\centering
\caption{Case III: The ML estimates of parameters $\theta$,  $\sigma$ and $\varepsilon$ from ESt distribution for $\varepsilon=-0.2, -0.5$ and $-0.8$.}
\label{simESt}
\scalebox{0.9}{
\begin{tabular}{cccccccc}\hline
& $\tau$ & $\hat{\tau}$ & ${Var}(\hat{\tau})$ & ${MSE}(\hat{\tau})$  & $\hat{\tau}$ & ${Var}(\hat{\tau})$ & ${MSE}(\hat{\tau})$  \\  \hline \hline
 &   & & $n=      30$ & &   &$n=      50$  &   \\ \hline\hline
$\theta$&$0.0$&$ -0.0631$&$  0.1751$&$  0.1779$&$ -0.0624$&$  0.1044$&$  0.1083$\\

$\sigma$&$1.0$&$  0.8609$&$  0.1210$&$  0.1395$&$  0.8768$&$  0.1071$&$  0.1223$\\

$\varepsilon$&$-0.2$&$ -0.2054$&$  0.0602$&$  0.0602$&$ -0.2316$&$  0.0364$&$  0.0374$\\
\hline\hline
 &   & & $n=     100$ & &  &$n=     150$  &   \\ \hline\hline
$\theta$&$0.0$&$  0.0021$&$  0.0617$&$  0.0617$&$ -0.0206$&$  0.0395$&$  0.0399$\\

$\sigma$&$1.0$&$  0.9573$&$  0.0391$&$  0.0411$&$  0.9867$&$  0.0266$&$  0.0268$\\

$\varepsilon$&$-0.2$&$ -0.2062$&$  0.0177$&$  0.0178$&$ -0.2147$&$  0.0134$&$  0.0137$\\
\hline \hline
 &   & & $n=      30$ & &   &$n=      50$  &   \\ \hline\hline
$\theta$&$0.0$&$  0.0270$&$  0.2015$&$  0.2023$&$ -0.0021$&$  0.0826$&$  0.0826$\\

$\sigma$&$1.0$&$  0.9368$&$  0.1295$&$  0.1335$&$  0.9169$&$  0.0808$&$  0.0889$\\

$\varepsilon$&$-0.5$&$ -0.5095$&$  0.0570$&$  0.0571$&$ -0.5081$&$  0.0325$&$  0.0326$\\ 
 \hline \hline
 &   & & $n=     100$ & &  &$n=     150$  &   \\ \hline\hline
$\theta$&$0.0$&$  0.0139$&$  0.0466$&$  0.0468$&$ -0.0199$&$  0.0361$&$  0.0365$\\

$\sigma$&$1.0$&$  0.9820$&$  0.0090$&$  0.0093$&$  0.9921$&$  0.0058$&$  0.0059$\\

$\varepsilon$&$-0.5$&$ -0.4874$&$  0.0138$&$  0.0140$&$ -0.5002$&$  0.0109$&$  0.0109$\\
\hline\hline
 &   & & $n=      30$ & &   &$n=      50$  &   \\ \hline\hline
$\theta$&$0.0$&$  0.1861$&$  0.1262$&$  0.1609$&$  0.0217$&$  0.0786$&$  0.0791$\\

$\sigma$&$1.0$&$  1.0403$&$  0.1442$&$  0.1458$&$  0.9810$&$  0.0510$&$  0.0514$\\

$\varepsilon$&$-0.8$&$ -0.7706$&$  0.0329$&$  0.0338$&$ -0.8153$&$  0.0205$&$  0.0207$\\
 \hline\hline
 &   & & $n=     100$ & &  &$n=     150$  &   \\ \hline\hline
$\theta$&$0.0$&$  0.0181$&$  0.0271$&$  0.0274$&$  0.0162$&$  0.0178$&$  0.0180$\\

$\sigma$&$1.0$&$  0.9481$&$  0.0468$&$  0.0495$&$  0.9956$&$  0.0265$&$  0.0266$\\

$\varepsilon$&$-0.8$&$ -0.8000$&$  0.0083$&$  0.0083$&$ -0.7984$&$  0.0052$&$  0.0052$\\
\hline\hline
\end{tabular}
}
\end{table}

Tables \ref{simESN}, \ref{simESL} and \ref{simESt} show that when the sample sizes increase the values of MSE decreases, as expected. In some cases, for examples, Case I with $\varepsilon=-0.5$, $n=150$ and Case II with $\varepsilon=-0.8$, $n=150$ for estimates of parameter $\theta$ show that the more data are needed to converge the value of true parameter. The same situations are observed for the estimates of parameters $\varepsilon$ in Case II with all values of skewness parameter for all sample sizes and $\sigma$ in Case III with $\varepsilon=-0.2$ and $-0.5$ for small sample sizes ($n=30$ and $n=50$) generally. The asymptotic variance values in tables 1-2 and the simulated MSE values are close to each other, as expected. From here, by using the well known properties of ML estimators, it is said that they are asymptotically efficient and asymptotic normally distributed (see theorem 4.1).

%
%
%
%
%
%

As a result, the simulation study can be preferable to see the precise performance of ML estimators, however using the KS test statistic and information criterions, such as AIC and BIC to see the fitting performance of the proposed functions is also needed. Due to this reason, we consider to give the KS test statistic, AIC and BIC which are inevitable tools to see the fitting performance of functions in the real data case. The subsection given below consists of examples from a real data sets.

\subsection{Modelling Real Data via Using Distributions}

\textbf{Examples}: The univariate cases of the data sets in cDNA micro array were modelled by \cite{Arslan09a,Arslan09b}. The parametric models are used to fit the data sets. This data sets were also analysed by \cite{Acitas13,PurdomHol05} for modelling them in the univariate case. In this study, we consider to model the cDNA micro array datasets. In these data sets, AT-matrix: 1416 genes $\times$ 118  drugs which are available at the web site

\begin{verbatim}http://discover.nci.nih.gov/
nature2000/data/selected_data/at_matrix.txt 
\end{verbatim}
are considered to test the performance of ESEP and its heavy tailed form (ESGT). Considering the data set in a class, such as micro array, is beneficial, because the data set in one class tried to be represented by means of PDFs in same class. The matching between the data set and the hypothetical models, that are a class of PDFs, can be given, in other words, the nature of phenomena must be represented by more appropriate hypothetical models. In this work, we tried to represent the epsilon skew and exponential power distribution that includes capacity for the heavy-tailed modelling as well. The variable named as "SID 416227, ESTs [5':, 3':W86124]" was fitted with these functions. Table 4 gives the estimates of parameters $\theta$, $\sigma$ and $\varepsilon$. To see fitting  capacity of models, the probability values of KS test statistic, represented by $P(KS)$, AIC and BIC from the models are computed and given in Table 6.
\begin{table}[!htb]
\label{esepfamilyone}
\begin{center}
\caption{Parameter estimates for the microarray data set via using the different parametric models}
\begin{tabular}{c c c c c c c c}
\hline Parametric Models &$\hat{\theta}$ & $\hat{\sigma}$  & $\hat{\varepsilon}$ & $P(KS)$ & $AIC$ & $BIC$ \\ [0.5ex]  
\hline \hline 
\multicolumn{7}{c}{Example 1} \\ \hline
ESL ($\alpha=1$) &  0.0321 &   0.0715  &  0.0687  & 0.5394 & -369.0365  &  -360.7244 \\
ESEP ($\alpha=1.01$)&  0.0317  &  0.0725 &   0.0691  &  0.4894 &  -366.8772  & -358.5651 \\
ESGT($\alpha=1,q=0.98$)  & 0.0270   & 0.0657 &  -0.0577  & 0.3134   &  -90.4600 & -82.1480 \\
 \hline  \hline
\multicolumn{7}{c}{Example 2} \\ \hline
 ESL ($\alpha=1$) &	 0.0120 &     0.0800   & 0.0284  & 0.6962 &  -342.6758 &  -334.3637\\
ESEP ($\alpha=1.01$) & 0.0120 &   0.0810 & 0.0287 & 0.6817 & -340.6565 & -332.3444 \\
ESGT($\alpha=1,q=0.95$)   &0.0130  &  0.0797 &   0.0002 &  0.1888 &  -52.8632 & -44.5512 \\
 \hline  \hline
\multicolumn{7}{c}{Example 3} \\ \hline
 ESL ($\alpha=1$) & 0.1059 &  0.0812&  0.3911   &  0.7007  &-339.0418  & -330.7297   \\
ESEP ($\alpha=1.01$) &   0.0440  & 0.0826 &  0.1205  &  0.6660 &  -335.8731 & -327.5610 \\
ESGT($\alpha=1,q=1.1$)   &  0.0339  & 0.0883 &  -0.000043   &  0.0212   &  -48.7423  & -40.4303 \\
 \hline  \hline
\multicolumn{7}{c}{Example 4} \\ \hline
 ESL ($\alpha=1$) &	 -0.0030  &  0.0752 &  -0.1228 & 0.8810  & -357.1366 &  -348.8245    \\
ESEP ($\alpha=0.9999$) &  0.0011  &  0.0766 &  -0.0618  &  0.6712 & -353.0269  &  -344.7149 \\
ESGT($\alpha=1,q=0.98$)   &   0.0060 &   0.0784  &  0.0027    &  0.0916   &  -63.1591  &   -54.8471 \\  
\hline
\end{tabular}
\label{sheeptable}
\end{center}
\end{table}


The following matrices show the Cram\'{e}r-Rao lower bound of the estimators for parameters $\theta$, $\sigma$ and $\varepsilon$ from the distributions. They are given in order from examples 1, 2 and 3. The variance-covariance matrix is for the estimators $\hat{\boldsymbol \tau}=(\hat{\theta},\hat{\sigma},\hat{\varepsilon})$ of parameters:
$$ Cov(\hat{\boldsymbol \tau}) = \begin{bmatrix}
  Var(\hat{\theta}) & Cov(\hat{\theta},\hat{\sigma}) & Cov(\hat{\theta},\hat{\varepsilon}) \\ & Var(\hat{\sigma}) & Cov(\hat{\sigma},\hat{\varepsilon})  \\
      &      &  Var(\hat{\varepsilon})  
\end{bmatrix}.$$

For the example 1, the following matrix is for ESEP with $\alpha=1.01$:
$$ Cov(\hat{\boldsymbol \tau}) = \begin{bmatrix}
     0.000178547813762        &           0  & -0.000878525941104 \\
       0  & 0.000044118987261  &            0  \\
 -0.000878525941104             &      0 &  0.008518791442917
\end{bmatrix}.$$

For the example 1, the following matrix is for ESGT($\alpha,q$) with $\alpha=1.00001$ and $q=0.98$:
$$ Cov(\hat{\boldsymbol \tau}) = \begin{bmatrix}
    0.000000155842497         &          0  & 0.000018048839655 \\
       0  & 0.000000000035809       &            0 \\
 0.000018048839655      &             0 &  0.008446385195929
\end{bmatrix}.$$

For the example 2, the following matrix is for ESEP with $\alpha=1.01$:
$$ Cov(\hat{\boldsymbol \tau}) = \begin{bmatrix}
   0.000223876064821      &             0 &  -0.000985691342202  \\
  0 &  0.000055100949914          &         0 \\
  -0.000985691342202    &               0 &  0.008552588217128
\end{bmatrix}.$$

For the example 2, the following matrix is for ESGT($\alpha,q$) with $\alpha=1.00001$ and $q=0.95$:
$$ Cov(\hat{\boldsymbol \tau}) = \begin{bmatrix}
   0.000000333448636           &        0 &  0.000026236717787  \\
      0  & 0.000000000171647        &           0 \\
 0.000026236717787          &         0 &  0.008474617300208
\end{bmatrix}.$$

For the example 3, the following matrix is for ESEP with $\alpha=1.01$:
$$ Cov(\hat{\boldsymbol \tau}) = \begin{bmatrix}
   0.000229927085055       &            0 & -0.000992045241996\\
                0  & 0.000057377675611    &              0  \\
  -0.000992045241996     &              0 &  0.008435215195203
\end{bmatrix}.$$

For the example 3, the following matrix is for ESGT($\alpha,q$) with $\alpha=1.00001$ and $q=1.1$:
$$ Cov(\hat{\boldsymbol \tau}) = \begin{bmatrix}
     0.000000540340673       &            0  & 0.000034631474732 \\
                   0 &  0.000000000353848         &          0  \\
     0.000034631474732          &         0  & 0.008474620646539
\end{bmatrix}.$$

For the example 4, the following matrix is for ESEP with $\alpha=0.9999$:
$$ Cov(\hat{\boldsymbol \tau}) = \begin{bmatrix}
   0.000197992173257       &            0  & -0.000914074903316\\
                   0  & 0.000049696265078       &            0\\
  -0.000914074903316          &         0  & 0.008441326498198     
\end{bmatrix}.$$

For the example 4, the following matrix is for ESGT($\alpha,q$) with $\alpha=1.00001$ and $q=0.98$:
$$ Cov(\hat{\boldsymbol \tau}) = \begin{bmatrix}
 0.000000316926028            &       0  & 0.000025781506306\\
       0 &  0.000000000147112  &                 0\\
         0.000025781506306         &          0 &  0.008474554487626
\end{bmatrix}.$$

Note that when $\alpha=1$, the results are not obtained for the ESGT($\alpha,q$) due to the function $\Gamma$. So, $\alpha$ was taken to be near to $1$.

The following figures are considered to illustrate the fitting performance of parametric models via the cumulative density function (CDF). It can be observed that the figures 1 and 2 include the outliers that are near to the underlying distribution, because the trend of end points that are left lower and upper sides of empirical CDF (blue line) in figures 1 and 2 is not high. However, one can observe that the trend of end points that are lower and upper sides of empirical CDF in the figure 3 is high. So, the data set in example 3 includes outlier(s) can not be observable from the empirical CDF. However, not only shape of empirical CDF but also the estimated value of skewness parameter from ESL distribution show that the data set of the example 3 includes asymmetry. The example 4 shows that ESL distribution has high fitting performance as it is seen from the value of  $P(KS)$. The upper side of empirical CDF for the example 4 shows there can be outliers in data set. In example 4, the asymmetry and the outliers that are small amount are observed together even if the potential outliers are near to the underlying distribution. In other side of our comment, the modelling capacity that is seen from high value of $P(KS)$ is an important indicator as well.

From the discussion given above, it is seen that instead of using histogram and superimposing the PDF on histogram in order to depict the asymmetry and the possible outlier(s) in data set, using CDF for discovering the nature of data can be more preferable in order to make a good matching between the data set and the parametric models as a reliable illustrating and also we can support this illustration via computing the value of $P(KS)$ to have the fitting competence of parametric model on data set. 
\begin{figure}[!htb] \label{fig:cdfexample}
\centering
    \includegraphics[width=.85\textwidth]{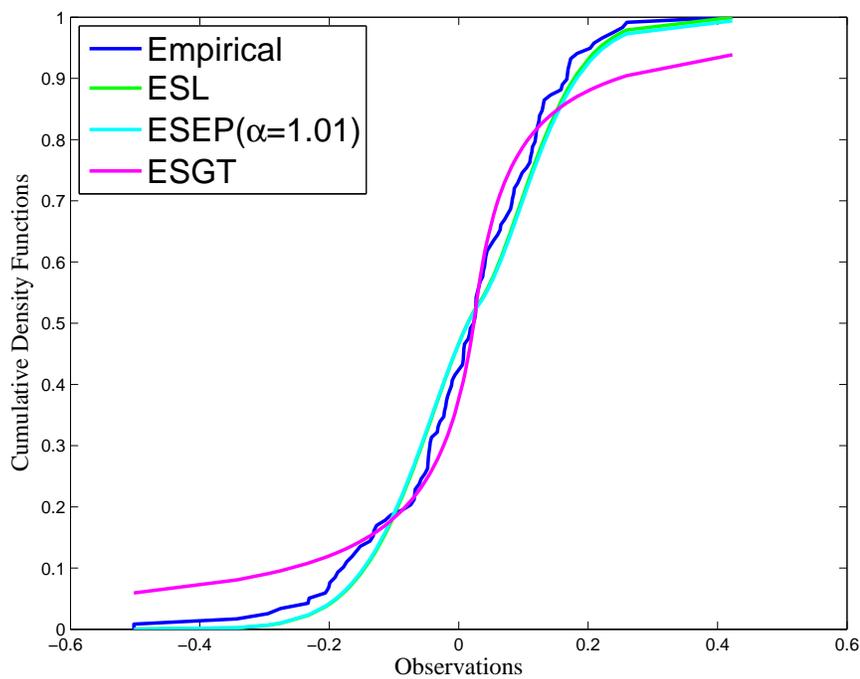}
 \caption{Empirical CDF of the data with the fitted CDF functions for microarray data of Example 1}
\end{figure}
\begin{figure}[!htb] \label{fig:cdfexample2}
\centering
    \includegraphics[width=.85\textwidth]{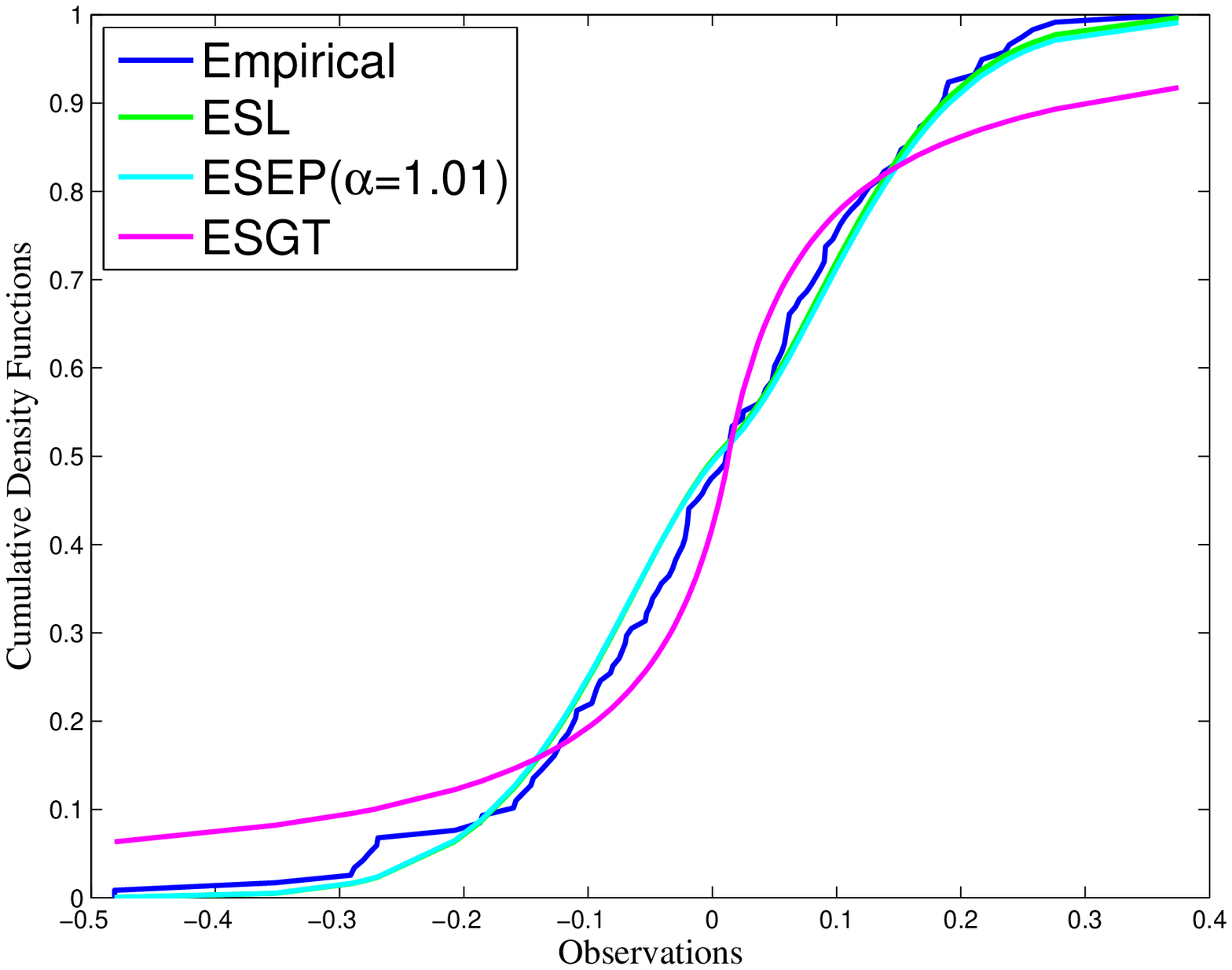}
 \caption{Empirical CDF of the data with the fitted CDF functions for microarray data of Example 2}
\end{figure}
\begin{figure}[!htb] \label{fig:cdfexample3}
\centering
    \includegraphics[width=.85\textwidth]{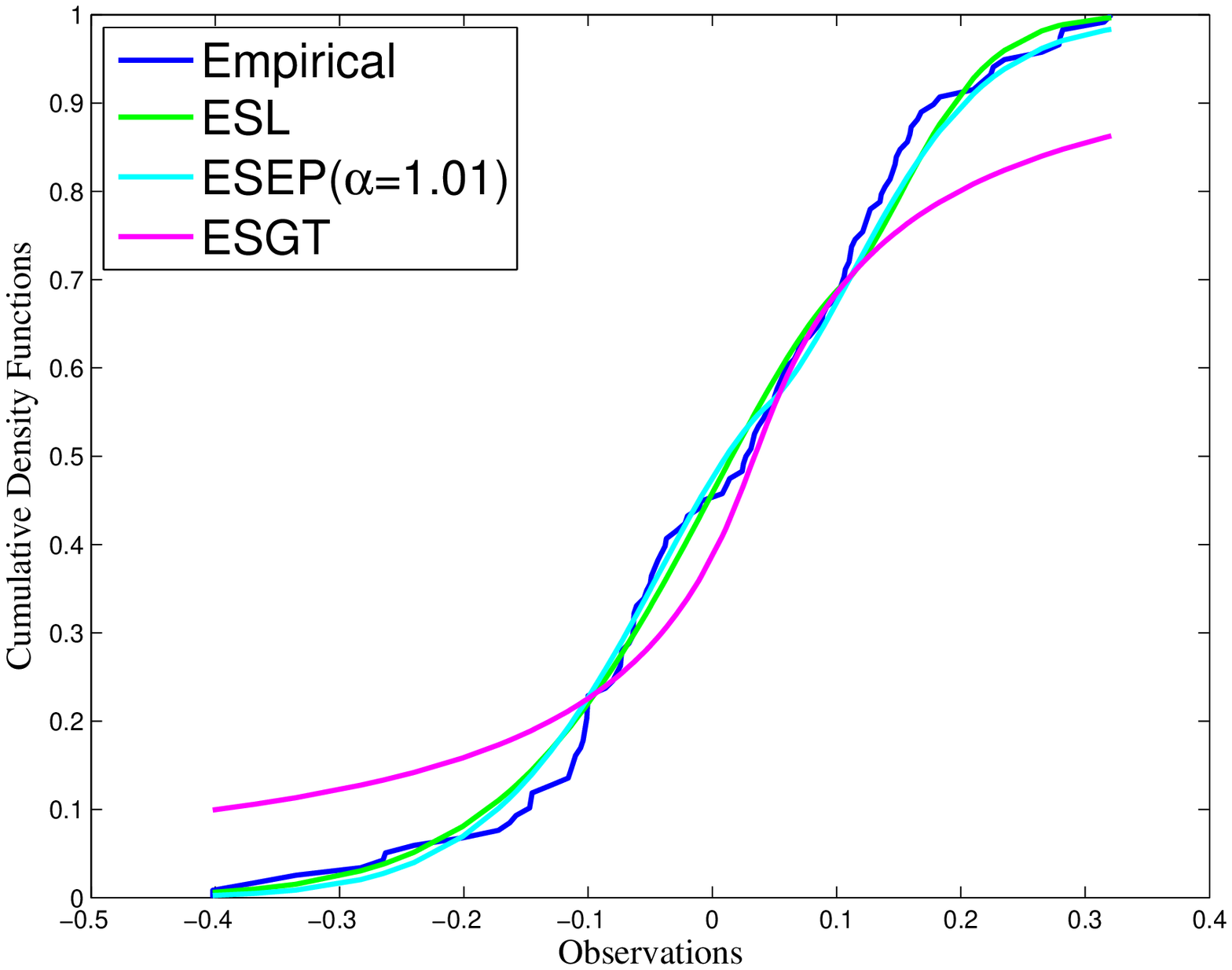}
 \caption{Empirical CDF of the data with the fitted CDF functions for microarray data of Example 3}
\end{figure}
\begin{figure}[!htb] \label{fig:cdfexample4}
\centering
    \includegraphics[width=.85\textwidth]{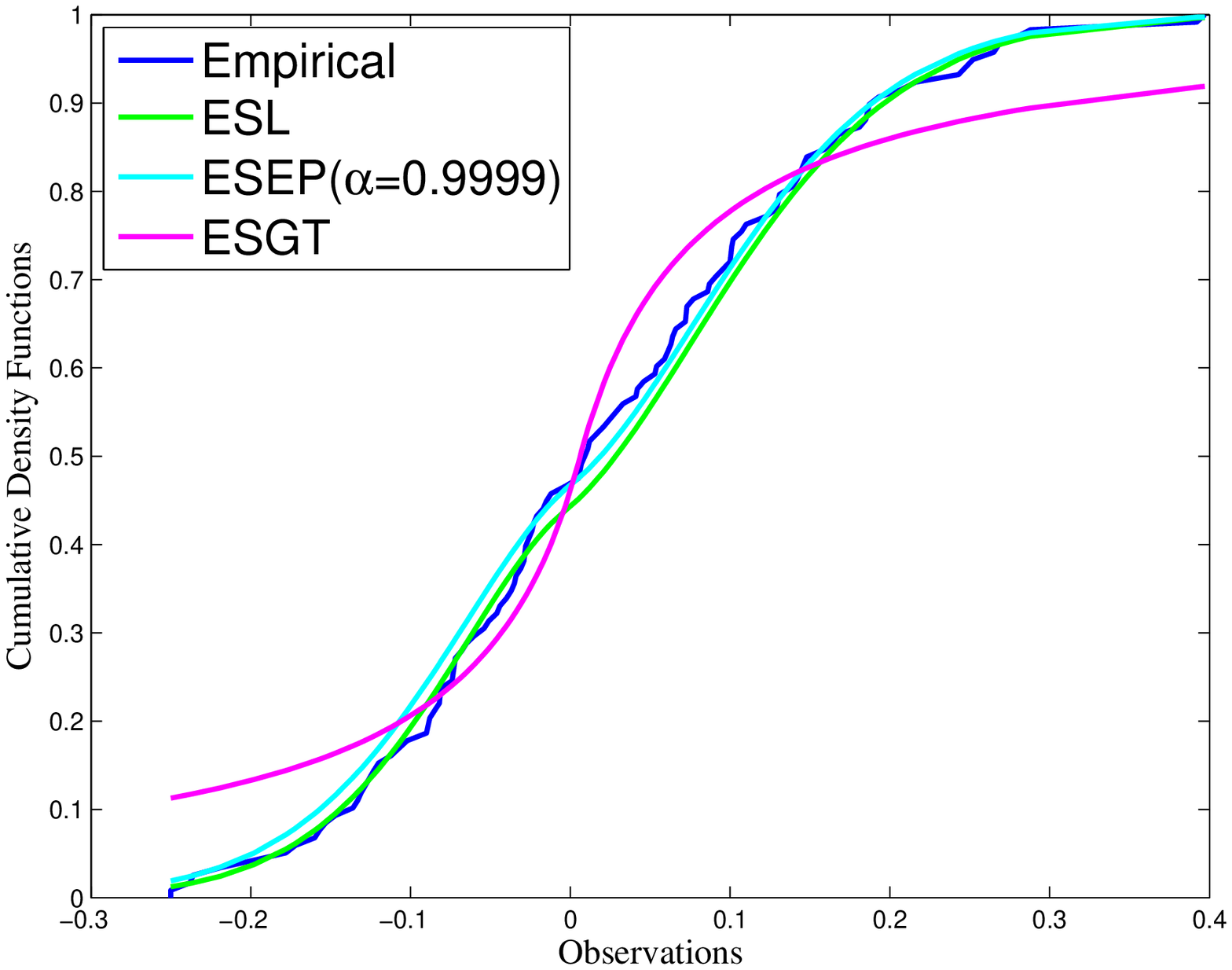}
 \caption{Empirical CDF of the data with the fitted CDF functions for microarray data of Example 4}
\end{figure}




\section{Conclusions}
The robustness properties of ML estimators of parameters location, scale and scale variant that is skewness  have been tested by means of the tools in Ref. \cite{Hampel86}. The ML estimators of these parameters in ESGT distribution have local robustness properties. However, the ML estimator of location parameter in ESEP with $\alpha \leq 1$ and ESGT distribution with $\alpha=2$ and $q=\nu/2$ has local and global robustness properties together. The ML estimators of location, scale and skewness parameters are given in the IRA form to compute the estimates of these parameters. The asymptotic properties of ML estimators for these parameters have been examined. We will examine the breakdown point of ML estimators of parameters $\sigma$ and $\varepsilon$ for arbitrary PDFs and also the robustness and asymptotic properties for linear and non-linear regressions. In future, we will prepare a R package to compute the estimates of parameters in these distributions that are univariate and multivariate forms.


\appendix
\section*{Appendix}
The PDF of ESEP is
\begin{verbatim}
f[\[Theta], \[Sigma], \[Epsilon], \[Alpha]] := 
 Log[\[Alpha]/(2^(3/2)*\[Sigma]*
      Gamma[1/\[Alpha]])] - ((Abs[
       x - \[Theta]]^\[Alpha])/(2^(\[Alpha]/2)*(1 - 
         Sign[x - \[Theta]]*\[Epsilon])^\[Alpha]*\[Sigma]^\[Alpha]))
\end{verbatim}
The Hessian matrix of parameters in ESEP distribution is
\begin{verbatim}
HessianESEP2 := 
 D[f[\[Theta], \[Sigma], \[Epsilon], \[Alpha]], {{\[Theta], \[Sigma], \
\[Epsilon], \[Alpha]}, 2}]
\end{verbatim}
The third order derivatives with respect to parameters are obtained by
\begin{verbatim}
HessianESEP3 := 
 D[HessianESEP2, {{\[Theta], \[Sigma], \
\[Epsilon], \[Alpha]}, 1}]
\end{verbatim}
The PDF of ESGT is
\begin{verbatim}
f[\[Theta], \[Sigma], \[Epsilon], \[Alpha], q] := 
 Log[(\[Alpha]*Gamma[1/\[Alpha] + q])/(Gamma[1/\[Alpha]]*
      Gamma[q]*2^(3/2)*q^(1/\[Alpha])*\[Sigma])] - (q + 1/\[Alpha])*
   Log[1 + ((Abs[
          x - \[Theta]]^\[Alpha])/(2^(\[Alpha]/2)*(1 - 
            Sign[x - \[Theta]]*\[Epsilon])^\[Alpha]*
         q*\[Sigma]^\[Alpha]))]
\end{verbatim}
The Hessian matrix of parameters in ESGT distribution is
\begin{verbatim}
HessianfESGT2 := 
 D[f[\[Theta], \[Sigma], \[Epsilon], \[Alpha], 
   q], {{\[Theta], \[Sigma], \[Epsilon], \[Alpha], q}, 2}]
\end{verbatim}
The third order derivatives with respect to parameters are obtained by
\begin{verbatim}
HessianfESGT3 := 
 D[HessianfESGT2, {{\[Theta], \[Sigma], \[Epsilon], \[Alpha], q}, 1}]
\end{verbatim}
One can obtain the Hessian matrices and third derivatives of log-likelihood functions w.r.t parameters for ESEP and ESGT distributions via using Mathematica 9.0.1.0 version.

\begin{verbatim}

\end{verbatim}

\end{document}